\documentclass{amsart}
\usepackage{graphicx}
\newtheorem{theorem}{Theorem}

\newtheorem{lemma}[theorem]{Lemma}

\newcommand{\abs}[1]{\left\vert#1\right\vert}

\usepackage{amsmath,amscd}
\begin{document}
\title[]{Totally Real Perturbations and Non-Degenerate Embeddings of $S^3$}%
\author{Ali M. Elgindi}%

\begin{abstract}
In this article, we demonstrate methods for the local removal and modification of complex tangents to embeddings of $S^3$ into $\mathbb{C}^3$. In particular, given any embedding of $S^3$ and a neighborhood of the complex tangents of the embedding, we show that there exists a ($\mathcal{C}^0$-close) totally real embedding which agrees with the original embedding outside the given neighborhood of the complex tangents. We also demonstrate that given any knot type $K \subset S^3$, either there exists an embedding of $S^3$ which assumes non-degenerate complex tangents exactly along $K$ or there exists a non-degenerate embedding complex tangent along two unlinked copies of $K$ (both cases may hold). We also note possible directions of future investigations.
\end{abstract}

\maketitle
\par\ \par\
\section{Introduction}
\par\ \par\
 Complex tangents to an embedding $M^k \hookrightarrow \mathbb{C}^n$ are points $x \in M$ so that the tangent space to $M$ at $x$ (considered as a subspace of the tangent space of $\mathbb{C}^n$) contains a complex line. If $k > n$, all points of $M$ are necessarily complex tangent, by virtue of the dimensions. If $k \leq n$, some (or all) points of $M$ may have strictly real tangent space. If all points of $M$ are real, we say that the embedding is totally real. In general, the dimension of the maximal complex tangent space of $M$ at $x$ is called the dimension of the complex tangent. An embedding is called CR (Cauchy-Riemann) if the complex dimension of all points of $M$ are the same. In [7], Gromov used the h-principle to demonstrate that the only spheres admitting totally real embeddings $S^n \hookrightarrow \mathbb{C}^n$ were $S^1$ and $S^3$. In dimension 3, CR-structures are either totally real or totally complex, two extremes which have been subject to much study. In his paper [6], Forstneric fully resolved the totally real problem in dimension 3; in particular he showed that every closed oriented 3-manifold can be embedded totally real into $\mathbb{C}^3$. 
\par\
The more general situations (proper subsets of complex tangents) have not been as thoroughly investigated. The work of Ahern and Rudin (in [1]) analyzed means of determining complex tangents for specific examples and Webster (in [12]) derived topological invariants for the set of complex tangents of embeddings of real $n$-manifolds into $n$-dimensional complex Euclidean space. In the general situation of embeddings of a real manifold into a complex manifold of arbitrary dimensions, the work of Lia in [9] and later of Domrin in [2] demonstrated formulas relating the characteristic classes of the set of complex tangents with those of the embedded manifold and the ambient complex manifold. Their results hold only under special assumptions on the dimensions of the manifolds and the structure of the complex tangents. Furthermore, Slapar in [11] considered the situation of a generic closed oriented manifold embedded into a complex manifold with codimension 2. In this situation, complex tangents generically are discrete (and finite), and can be classified as either elliptic or hyperbolic by comparing the orientation of the tangent space of the embedded manifold at the complex tangent (which is necessarily a complex vector space) with the induced orientation of the tangent space as a complex subspace of the tangent space of the complex ambient manifold. He demonstrated that every pair of complex tangents consisting of one elliptic point and one hyperbolic point can be canceled out; in particular there exists an arbitrarily small isotopy of embeddings so that the set of complex tangents of the isotoped embedding consists exactly of the complex tangents of the original embedding minus the given pair of complex tangents (one elliptic and one hyperbolic). 
\par\
In this article, we wish to investigate the means (and flexibility) of deforming embeddings of $S^3 \hookrightarrow \mathbb{C}^3$ to locally remove or manipulate the set complex tangents. In this setting, complex tangents generically arise as knots (links). In our paper [3] we proved that for every knot type in $S^3$ there exists an embedding $S^3 \hookrightarrow \mathbb{C}^3$ which assumes complex tangents exactly along a knot of the prescribed type, with one point being degenerate. 
We further demonstrated that for every knot type, there exists an embedding whose complex tangents form a wholly non-degenerate knot of the prescribed type union an isolated degenerate point. Our work could be readily extended to include all types of links. 
\par\
In our paper [5], we demonstrated an obstruction to (locally) removing an isolated complex tangent to an embedding $M^3 \hookrightarrow \mathbb{C}^3$, leaving the embedding unchanged outside a small neighborhood of the degenerate point (and without adding any new complex tangents). The obstruction is found to be the homotopy class of the image of the  
boundary 2-sphere $S$ of a sufficiently small neighborhood of the degenerate point under the Gauss map in $\pi_2(\mathbb{Y})$, where $\mathbb{Y}$ is the submanifold of the real Grassmannian $\textbf{\emph{G}}_{3,6}$ consisting of totally real 3-planes in $\mathbb{R}^6 = \mathbb{C}^3$. Using the h-principle, we demonstrated that the isolated complex tangent can be (locally) removed precisely when this obstruction vanishes. Further, (with the vanishing of the obstruction) the resulting embedding can be taken to be $\mathcal{C}^0$-close to the original embedding (in the perturbed neighborhood). 
\par\
We will use our findings from [5] to obtain two results for the case $M = S^3$. Our first result extends the previous work of Gromov and Forstneric (among others) to find that given any embedding of $S^3$ with (not dense) complex tangent set $\mathit{C} \subset S^3$ and an open ball $\mathit{U} \subset S^3$ containing $\mathit{C}$, there exists a ($\mathcal{C}^0$-close) totally real embedding of $S^3$ that agrees with the original given embedding outside of $\mathit{U}$.
\par\
For our second result we demonstrate  that given any knot type $K$ there exists an embedding $S^3 \hookrightarrow \mathbb{C}^3$ which either assumes non-degenerate complex tangents precisely upon $K$ or along two identical copies of the knot $K$ (unlinked). The complex tangents can be taken to be wholly non-degenerate in either case (we note that both cases may hold, at least for some knot types). This is an extension of our result in [3] where we constructed embeddings assuming complex tangents non-degenerate along a knot of any given type, union an isolated (degenerate) point. This result will also hold for $K$ being a link (type), by an extension of our work in [3].
\par\
\section{On Totally Real Perturbations}
\par\
In our first consideration of interest, let $S^3 \hookrightarrow \mathbb{C}^3$ be an embedding of class $\mathcal{C}^k$ and let $\mathit{C} \subset S^3$ be the set of complex tangents. In the generic situation, $C$ will be of codimension 2 in $S^3$, and as such will in the form of a knot (or link). But here it is not necessary to assume the embedding is generic; let us assume only that there exists an open ball $\mathbf{U}$ containing $\mathit{C}$ (some points are complex tangent, but are not dense). Let $\mathbf{S}$ be the boundary of $\mathbf{U}$, in particular $\mathbf{S} \subset S^3$ will be a 2-sphere. 
\par\
We recall from our work in [5]:
\begin{theorem}{([5])}
 Suppose an embedding of a 3-manifold $f: M \hookrightarrow \mathbb{C}^3$ (of class  $\mathcal{C}^k$) with Gauss map $G$ admits an isolated complex tangent $x \in M$ for which there is a neighborhood $\textbf{B}$ of $x$ whose boundary sphere $\textbf{S} \subset M$ satisfies: $[G(\textbf{S})] = 0 \in \pi_2(\mathbb{Y})$. Then there exists a $\mathcal{C}^k$-embedding $\widetilde{f}: M \hookrightarrow \mathbb{C}^3$ that can be taken to be $\mathcal{C}^0$-close to $f$ and regularly homotopic to $f$ so that:
 \par\
 1. $f = \widetilde{f}$ on $M \setminus \textbf{B}$
 \par\
 2. The set of complex tangents of $\widetilde{f}$ equals the set of complex tangents of $f$ minus the point $x$, i.e.:
 $\aleph_{\widetilde{f}} = \aleph_{f} \setminus \{x\}$
      \end{theorem}
 It is direct to see that the proof of the above theorem (in [5]) extends analogously when we replace an isolated point by a (bounded) set, in particular we can show that we may "fill in" the neighborhood $\mathbf{U}$ of the set of complex tangents $C$ in such a totally real manner if and only if $[G(\mathbf{S})] = 0$ in $\mathbb{Y}$. More precisely, with the vanishing of the obstruction, there exists an embedding of $S^3$ (of class $\mathcal{C}^k$) which is totally real in $\mathbf{U}$ and which agrees with the original embedding in the complement of $\mathbf{U}$, without adding any new complex tangents. 
 \par\
By construction, the original embedding is totally real in the complement of $\mathbf{U}$. Hence, we found an embedding of $S^3$ which is totally real, and which agrees with the original given embedding outside of the ball $\mathbf{U}$ (and $\mathcal{C}^0$-close throughout). However, this only holds with the assumption that $[G(\mathbf{S})] = 0$ in $\mathbb{Y}$, where $\mathbf{S}$ is the boundary of $\mathbf{U}$ and $G$ is the Gauss map of the embedding. This is in fact true for all embeddings of $S^3$ in complex Euclidean space.
\begin{lemma}: Let $S^3 \hookrightarrow \mathbb{C}^3$ be an embedding and let $\mathbf{U}$ be a 3-ball in $S^3$ containing all the complex tangents of the embedding, and consider its boundary sphere $\mathbf{S}$. Then the image of the Gauss map restricted to $\mathbf{S}$ is trivial in the homotopy group $\pi_2 (\mathbb{Y})$.
      \end{lemma}
\par\
$\textbf{\emph{Proof:}}$
\par\
As $\mathbf{S} \subset S^3$ is a 2-sphere, it will necessarily cut $S^3$ into two parts, the "internal" 3-ball which is the given $\mathbf{U}$, and the "external" 3-ball $\mathbf{E}$, which is the interior of the complement of $\mathbf{U}$ in $S^3$.
The 2-sphere $\mathbf{S}$ can then be contracted to a point in the external 3-ball $\mathbf{E}$, in particular there exists a homotopy $H:S^2 \times I \rightarrow S^3$ with $Im(H_t) \subset \mathbf{E}$ ($t>0$) so that:
\par\
$H|_0 (S^2) = \mathbf{S}$ and $H|_1 (S^2) = *$, where $*$ is a fixed point in the 3-ball $\mathbf{E}$. Consider now the homotopy: $K = G \circ H:S^2 \rightarrow \mathbb{Y}$.
\par\
Note that the image of the homotopy $K$ is indeed contained in $\mathbb{Y}$ as the image of $H$ is contained in the external 3-ball $\mathbf{E}$ , which contains no complex tangents of the embedding.
\par\
Now, $K|_0 (S^2)= G(\mathbf{S})$ and $K|_1 (S^2) = G(*) = p \in \mathbb{Y}$, a point.
\par\
Hence, we find that $[G(\mathbf{S})] = [p] = 0$ in $\pi_2(\mathbb{Y})$, which is what was claimed.
\par\                             
$\textbf{\emph{QED}}$
\par\ \par\
Hence, by our work in [5] and our arguments above preceding the Lemma, we have proven:
\begin{theorem} Let $F: S^3 \hookrightarrow \mathbb{C}^3$ be a $\mathcal{C}^k$-embedding and let $\mathbf{U}$ be a 3-ball in $S^3$ containing all the complex tangents of the embedding. Then there exists  a totally real embedding $\widetilde{F}: S^3 \hookrightarrow \mathbb{C}^3$ of class $\mathcal{C}^k$ so that $\widetilde{F}$ is $\mathcal{C}^0$-close to $F$ over $\mathbf{U}$ and $\widetilde{F} = F$ over $S^3 \setminus \mathbf{U}$.
      \end{theorem}
      \par\
We note our result is a generalization of the works of Gromov and Forstneric, who asserted the existence of totally real embeddings  of $S^3$ (and other closed 3-manifolds). Here, we have demonstrated that any embedding of $S^3$ which is not "densely complex" may be $\mathcal{C}^0$-closely perturbed to be totally real.
\par\
A possible future question would be to ask if all embeddings (even those which are complex) can be closely perturbed to be totally real. We may even ask this question for general manifolds, although we note that our arguments above for $S^3$ could not possibly generalize to other 3-manifolds, as the 3-sphere is the only closed 3-manifold which may be split into two balls by a 2-sphere (see Alexander's Theorem in [8]).
\par\
\section{Knots of Complex Tangents (and their doubles)}
\par\
We now recall our result from [3], where we demonstrated that for any given knot type (represented by $K \subset S^3$) and $k \in \mathbb{N}$, there exists a  graphical embedding $f:S^3 \rightarrow \mathbb{C}$ of class $\mathcal{C}^k$ which assumes non-degenerate complex tangents precisely along the knot $K$ union an isolated (degenerate) point $*$. 
\par\
While we would ultimately hope to be able to remove the degenerate point while leaving the knot $K$ unaffected, we are not able at this time to fully guarantee the existence of such an embedding with complex tangents wholly non-degenerate exactly along (a knot of type) $K$. Instead, we are able to demonstrate that if there exists no such embedding of $S^3$ with a knot of the prescribed type alone as complex tangent (non-degenerate), then we can make a complex "double" of the knot, in the sense that we can construct an embedding assuming complex tangents precisely along two identical (unlinking) copies of a knot of the given type, which are again non-degenerate.
\par\
More precisely:
\begin{theorem} Consider any topological knot type in $S^3$ and let $k \in \mathbb{N}$. Then (at least) one of the following two cases holds:
\par\
1. There exists a $\mathcal{C}^k$-embedding $S^3 \hookrightarrow \mathbb{C}^3$ assuming complex tangents precisely along a knot of the prescribed type which is wholly non-degenerate.
\par\
2. There exists a $\mathcal{C}^k$-embedding $S^3 \hookrightarrow \mathbb{C}^3$ assuming complex tangents precisely along two identical copies of a knot of the prescribed type (one in the upper hemisphere and one in the lower hemisphere), which are hence unlinked and are wholly non-degenerate.
      \end{theorem}
\par\
$\textbf{\emph{Proof:}}$
\par\
Consider a knot type in $S^3$. We constructed a $\mathcal{C}^k$-function $f:\mathbb{C}^2 \rightarrow \mathbb{C}$ in our article [3] whose graph (restricted to $S^3$)       
 $F:S^3 \hookrightarrow \mathbb{C}^3$ assumes complex tangents precisely along a non-degenerate knot of the prescribed type union an isolated degenerate point.
 We recall that the function $f$ is obtained by: 
 \par\
 $f(z,w)=(g \circ \psi)(z,w) \cdot (1-w)^r$
 \par\
 where $g:\mathbb{C}^2 \rightarrow \mathbb{C}$ is a (non-holomorphic) polynomial whose graph over the Heisenberg group $\mathbb{H}$ has complex tangents exactly along a knot $\widetilde{K}$ and $\psi: S^3 \setminus \{(0,1)\} \rightarrow \mathbb{H}$ is the standard biholomorphism. In particular, denoting $\varphi = \psi^{-1}$, the knot $\widetilde{K}$ is taken so that $K = \varphi(\widetilde{K})$ is a knot of our desired type in $S^3$, and in fact $\widetilde{K}$ can be assumed (without loss) to be wholly non-degenerate as the set of complex tangents to $graph(g)$ in $\mathbb{H}$. This follows from the fact that the tangential CR-operator to $\mathbb{H}$ is onto, which allows for flexibility to manipulate the complex tangents of embeddings $\mathbb{H} \hookrightarrow \mathbb{C}^3$ directly (in particular, via small perturbations). The function $f$ will then have graph $F$ over $S^3$ being a $\mathcal{C}^k$-embedding with complex tangents along $K$ (which will be non-degenerate) union an isolated (necessarily degenerate) complex tangent point at $(0,1)$, by the formulation of $f$.
\par\
(We refer the reader to our article [3] for the details of our construction)
\par\ \par\
We may further assume without loss of generality that the knot $K$ is contained in the lower-half hemisphere $D_- = \{(z,w) \in S^3 | Im(w) < 0$\}, again due to the flexibility to first manipulate the complex tangents in $\mathbb{H}$ using a scaling factor. Let $D^+, D^-$ denote the upper and lower hemispheres of $S^3$, respectively.
\par\ \par\
We may now consider the Gauss map of the embedding $F$, which we denote: $G: S^3 \rightarrow \textbf{G}_{3,6}$. By construction, this embedding will be non-degenerate along the knot $K$ and degenerate at $(0,1)$. Hence, the Gauss map will intersect the subspace $\mathbb{W} \subset \textbf{G}_{3,6}$ transversely along $K$ and intersect at $(0,1)$ in a degenerate manner (see [12]).
\par\ \par\
Now, let $\widetilde{\psi}: S^3 \setminus \{(0,-1)\} \rightarrow \mathbb{H}$ be the "inverted" biholomorphism, namely $\widetilde{\psi} = \psi \circ R$, where $R(z,w)=(z,-w)$. Consider then the embedding $\widetilde{F}: S^3 \hookrightarrow \mathbb{C}$ given by the graph of the function: 
\par\
$\widetilde{f} = (g \circ \widetilde{\psi}) \cdot (1+w)^r$.
\par\
By an elementary computation, we find that the complex tangents to $\widetilde{F} = graph(\widetilde{f}|_{S^3})$ will be precisely $K^* \cup \{(0,-1)\}$, where $K^* = \{(z, -w) | (z,w) \in K\}$, a reflection of the knot $K$ in the upper-hemisphere of $S^3$. Furthermore, the complex tangents will be non-degenerate along $K^*$ and for the same choice of $r$ as above in the construction of $f$, the embedding $\widetilde{F}$ will be of class $\mathcal{C}^k$.
\par\ \par\
Let $\epsilon > 0$ be so that $K \subset O_{2\epsilon} = \{(z,w) \in S^3 | Im(w) < -2\epsilon\}$. Then it is clear that $K^* \subset O^{2\epsilon} = \{(z,w) \in S^3 | Im(w)  > 2\epsilon\}$. Consider now the 2-sphere $S_\epsilon$ which is the boundary of the 3-ball $O_\epsilon$ and $S^\epsilon$ which is the boundary of $O^\epsilon$.
\par\ \par\
Now, the embedding $F: D^- \hookrightarrow \mathbb{C}^3$ will assume (non-degenerate) complex tangents exactly along $K$, while $\widetilde{F}: D^+ \hookrightarrow  \mathbb{C}^3$ will assume (non-degenerate) complex tangents exactly along $K^*$. Let $G, \widetilde{G}$ be the corresponding Gauss maps. 
\par\ \par\
We will then have two classes $[G(S_\epsilon)], [\widetilde{G}(S^\epsilon)]$ in the homotopy group $\pi_2(\mathbb{Y}) = \mathbb{Z}_2$. If either of these classes is trivial ($=0$), then by our work in [5] we may extend the corresponding embedding to the remaining hemisphere in a totally real fashion to obtain an $\mathcal{C}^k$-embedding $S^3 \hookrightarrow \mathbb{C}^3$ assuming (non-degenerate) complex tangents exactly along the knot $K$ (or the equivalent $K^*$) of the prescribed type. This is the first case of the theorem.
\par\ \par\
Let us now proceed with the second case of the theorem. Assume then that the above assumption does not hold, namely that neither of the classes $[G(S_\epsilon)], [\widetilde{G}(S^\epsilon)]$ are trivial in $\pi_2(\mathbb{Y}) = \mathbb{Z}_2$. Then they must necessarily be equal to each other (as the homotopy group has only two elements). As such, the maps $G: S_\epsilon \rightarrow \mathbb{Y}$ and $\widetilde{G}: S^\epsilon \rightarrow \mathbb{Y}$ are homotopic. We may write the homotopy as: $H: S^2 \times (-\epsilon,\epsilon) \rightarrow S_{1-\epsilon^2}^2 \times (-\epsilon,\epsilon) \rightarrow \mathbb{Y}$, with $S_{1-\epsilon^2}^2$ being envisioned as the 2-sphere of radius $1-\epsilon^2$. We have that: $H_{-\epsilon} = G|_{S_\epsilon}, H_\epsilon = \widetilde{G}|_{S^\epsilon}$.
\par\ \par\
Let $D_{2\epsilon} = \{(z,w) \in S^3 | -2\epsilon < Im(w) < 2\epsilon\}$. Consider the map $\mathcal{H}: D_{2\epsilon} \rightarrow \mathbb{Y}$ defined by: 
\par\
$\mathcal{H} (z,w) = G(z,w)$ for $-2\epsilon < Im(w) \leq -\epsilon$
\par\
$\mathcal{H} (z,w) = H_t (z,w), t=Im(w)$, for $-\epsilon < Im(w) <\epsilon$
\par\
$\mathcal{H} (z,w) = \widetilde{G} (z,w)$, for $\epsilon \leq Im(w) < 2\epsilon$.
\par\
Hence, $\mathcal{H}: D_{2\epsilon} \rightarrow \mathbb{Y}$ is a continuous map from a "thickened" 2-sphere into $\mathbb{Y}$. 
\par\ \par\
Let us now consider $D_{\epsilon}$, which we may write as: $D_{\epsilon}= S^2 \times I$, for the interval $I = [-\epsilon,\epsilon]$.
\par\ \par\
Consider the set of homotopy classes of maps $[S^2 \times I,\mathbb{Y};G|_{S_\epsilon} \cup \widetilde{G}|_{S^\epsilon}]$ that restrict along the boundary to the map $\mathcal{H}$ on $S^2 \times \{-\epsilon, \epsilon\}$, which is given by $G|_{S_\epsilon}, \widetilde{G}|_{S^\epsilon}$. As this set is non-empty (we constructed above the map $\mathcal{H}$), the set is in 1-1 correspondence with $\pi_3 (\mathbb{Y})$, as $(S^2 \times I) / \partial = S^3$ (see Hatcher in [8]). 
\par\ \par\  
Consider now the new embedding: $E: S^3 \hookrightarrow \mathbb{C}^3$, given by:
\par\
$F(z,w) = graph(f(z,w))$, if $Im(w) \leq -\frac{\epsilon}{2}$
\par\
$\widetilde{F} (z,w) = graph(\widetilde{f}(z,w))$, if $Im(w) \geq \frac{\epsilon}{2}$
\par\
$graph(\frac{1}{\epsilon} ((\frac{\epsilon}{2} - Imw)f(z,w) + (\frac{\epsilon}{2} + Imw) \widetilde{f}(z,w)))$, if $-\frac{\epsilon}{2} < Imw < \frac{\epsilon}{2}$.
\par\ \par\
We notice that since the embedding $E$ corresponds with $F$ in a neighborhood of $O_\epsilon$, the Gauss map of $E$, which we denote by $G^E$, has the property that: 
\par\
$G^E |_{S_\epsilon} = G|_{S_\epsilon}$.
\par\
In direct analogy (in upper-hemisphere) $G^E |_{S^\epsilon} = \widetilde{G}|_{S^\epsilon}$.
\par\ \par\   
Our goal is to demonstrate that the Gauss map $G^E:S^2 \times I \rightarrow \textbf{\textit{G}}_{3,6}$ is homotopic to a map $\widehat{G}:S^2 \times I \rightarrow \mathbb{Y}$ that agrees with $G^E$ on the boundary spheres $S_\epsilon, S^\epsilon$, and through such maps (agreeing on $S_\epsilon \cup S^\epsilon$). 
\par\
Note that $G^E \in [S^2 \times I, \textbf{\textit{G}}_{3,6};G|_{S_\epsilon} \cup \widetilde{G}|_{S^\epsilon}] \cong \pi_3(\textbf{\textit{G}}_{3,6})$ (see Hatcher in [8])
\par\ \par\
Now, let $V_{3,6}$ be the generalized Stiefel manifold of all 3-frames in $\mathbb{R}^6$ and consider the subset of totally real 3-frames ${V_{3,6}}^{tr}$. We get the natural maps sending a frame to the subspace spanned by the frames in $\mathbb{R}^6$:
\par\
$s: V_{3,6} \rightarrow \textbf{\textit{G}}_{3,6}$
\par\
$r: {V_{3,6}}^{tr} \rightarrow \mathbb{Y}$ 
\par\
(see Forstneric in [6] for reference for our work here and below)
\par\ \par\
Notice then the commutative diagram:
\par\ \par\ 
  $\begin{CD}
{V_{3,6}}^{tr} @>>{i_V}>            V_{3,6} \\
 @VVrV @VVsV \\
\mathbb{Y}  @>>{i_G}>          \textbf{\textit{G}}_{3,6}
\end{CD}$
\par\ \par\
where $i_V,i_G$ are inclusion maps.
\par\ \par\
Applying $\pi_3$, we get the resulting commutative diagram of groups:
\par\ \par\
 $\begin{CD}
\pi_3({V_{3,6}}^{tr}) @>>{{i_V}}_*>            \pi_3(V_{3,6}) \\
 @VVr_*V @VVs_*V \\
\pi_3(\mathbb{Y})  @>>{{i_G}_*}>          \pi_3(\textbf{\textit{G}}_{3,6})
\end{CD}$
\par\ \par\
Now, it is known that $\pi_3(V_{3,6}) = \mathbb{Z}_2$, and from the long exact sequence of the fiber bundle over $\textbf{\textit{G}}_{3,6}$, we get that the map $s_*$ is an isomorphism.
\par\ \par\
Also, notice that ${V_{3,6}}^{tr} \cong Gl_3(\mathbb{C})$ retracts to $U(3)$, and that the space of orthonormal 3-frames ${V_{3,6}}^O$ is a retract of $V_{3,6}$. Hence, the map $i_V: {V_{3,6}}^{tr} \rightarrow V_{3,6}$ is homotopy equivalent to the map $\alpha: U(3) \rightarrow {V_{3,6}}^O$, which is the composition of the inclusion $U(3) \hookrightarrow O(6)$ with the quotient map $O(6) \rightarrow {V_{3,6}}^O$. Since the first map is an isomorphism on $\pi_3$ and the second is onto on $\pi_3$ (see [7]), we see that $\alpha_*:\pi_3(U(3)) \rightarrow \pi_3({V_{3,6}}^O)$ is onto, and hence the inclusion map $(i_V)_*$ is an onto map on $\pi_3$.
\par\ \par\
Hence, as the above diagram in $\pi_3$ commutes, $s_* \circ {i_V}_* = {i_G}_* \circ r_*$ is onto, and since $\pi_3(Y) = \mathbb{Z}_2 = \pi_3(\textbf{\textit{G}}_{3,6})$ (as shown above), we have that:
\par\
${i_G}_*: \pi_3(\mathbb{Y})  \rightarrow \pi_3(\textbf{\textit{G}}_{3,6})$ is an isomorphism.
\par\
Recall that: 
\par\
$\pi_3(\textbf{\textit{G}}_{3,6}) \cong [S^2 \times I, \textbf{\textit{G}}_{3,6};G|_{S_\epsilon} \cup \widetilde{G}|_{S^\epsilon}]$ and $\pi_3(\mathbb{Y}) \cong [S^2 \times I, \mathbb{Y};G|_{S_\epsilon} \cup \widetilde{G}|_{S^\epsilon}]$.
\par\ \par\
Now, the Gauss map $G^E$ of the embedding $E$ is clearly an element of $[S^2 \times I, \textbf{\textit{G}}_{3,6};G|_{S_\epsilon} \cup \widetilde{G}|_{S^\epsilon}]$, and by the above equivalence and the fact that $i_G$ is the inclusion map, there exists a map $\widehat{G} \in [S^2 \times I, \mathbb{Y};G|_{S_\epsilon} \cup \widetilde{G}|_{S^\epsilon}]$ which is homotopic to $G$ through maps $G_t$ which all agree on the boundary $S^2 \times \{-\epsilon,  \epsilon\}$
\par\ \par\
Consider then the maps: $\mathcal{G}_t: S^3 \rightarrow \textbf{\textit{G}}_{3,6}$, where ${\mathcal{G}_t}|_{S^2 \times I} = G_t,  {\mathcal{G}_t}|_{S^3 \setminus (S^2 \times I)} = G^E$. Then ${\mathcal{G}}_0 = G^E$ and ${\mathcal{G}_1}|_{S^2 \times I} = \widehat{G}$, which is totally real on $S^2 \times I$. All these maps can be taken to be continuous. We may further extend all these maps continuously to the "thickened" annulus $D_{2\epsilon}$ using the original Gauss maps $G, \widetilde{G}$ on the lower and upper part, respectively.
\par\ \par\
In fact, by the Whitney Approximation Theorem, we may assume that all these maps are of class $\mathcal{C}^k$ (by taking a small homotopy). Furthermore, we may make use the relative version of the Whitney Approximation Theorem  and assume without loss of generality that $\mathcal{G}_t$ agrees with the  original Gauss maps $G, \widetilde{G}$ on the bottom and top parts (respectively) of $\textbf{A}_{2\epsilon} = D_{2\epsilon} \setminus D_{\epsilon}$. (See Lee in [10] for an exposition on the Whitney Approximation Theorem).
\par\ \par\
Hence, we have constructed a formal totally real map, which we now denote $\widehat{G} = \mathcal{G}_1: D_{2\epsilon} \rightarrow \mathbb{Y}$ of class  $\mathcal{C}^k$ that agrees with the Gauss maps $G, \widetilde{G}$ of the original embeddings on the bottom and top parts (respectively) of the neighborhood $\textbf{A}_{2\epsilon}$ of the boundary, and which is homotopic to $G^E$ through maps $D_{2\epsilon} \rightarrow \textbf{\textit{G}}_{3,6}$, In the language of the literature, we have shown that our constructed embedding $E$ is $\textit{formally totally real}$ (on $D_{2\epsilon}$).
\par\ \par\
Now, applying the h-principle for extensions to this map $E$ on $D_{2\epsilon}$, we see that since this is a formal totally real embedding which is holonomic for  $\frac {3\epsilon}{2} \leq \abs{Im(w)} < 2\epsilon$, there exists a holonomic solution $F^@: D_{2\epsilon} \hookrightarrow \mathbb{C}^3$ which agrees with the original $F$ for $Im(w) \geq \frac {3\epsilon}{2}$ and with $\widetilde{F}$ on $Im(w) \leq - \frac {3\epsilon}{2}$.
Furthermore, $F^@: D_{2\epsilon} \hookrightarrow \mathbb{C}^3$ is totally real and can be taken to be $\mathcal{C}^0$-close to the (graphical) embedding $E$.
\par\
 (see our work in [5] and Gromov in [7] for reference)
\par\ \par\
We may extend $F^@ : S^3 \hookrightarrow \mathbb{C}^3$ to the entire of $S^3$ in the natural way, namely: 
\par\
$F^@ = F$ for $Im(w) \leq -2\epsilon$
\par\
$F^@ = \widetilde{F}$ for $Im(w) \geq 2\epsilon$
\par\
And as (consistently) defined above on $D_{2\epsilon}$.
\par\ \par\
The fact that $F^@$ is indeed an embedding of $S^3$ can be seen implicitly as $F^@$ can be taken to be $\mathcal{C}^0$ close to $E$ on the lower and upper caps of $S^3$ ($Imw > 2\epsilon, Imw < -2\epsilon)$. It then follows that there can be no self-intersections of $F^@$ as it can be taken to be arbitrarily close to a (fixed) graphical embedding.
\par\ 
Furthermore, by construction the $\mathcal{C}^k$-embedding $F^@ : S^3 \hookrightarrow \mathbb{C}^3$ takes complex tangents precisely along the unlinked knots $K \cup K^*$, which are both non-degenerate.
\par\ \par\
As this is based on the assumption of the second case, we have demonstrated the results of the theorem.							
\par\
$\textbf{\emph{QED}}$
\par\ \par\
We note that the results of the above Theorem may be readily extended to include all classes of links, by an immediate generalization of our work in [3] (given that "all links are algebraic").
\par\ \par\
We further note that there are examples in which both cases of the Theorem will hold. In particular, in  our paper [4], we demonstrated two examples of embeddings of $S^3$, one assuming complex tangents along a circle and another assuming complex tangents along two copies of a circle (unlinked). In both examples the complex tangents were all non-degenerate.
\par\ \par\
To establish the existence of both cases, in particular for the existence of embeddings with non-degenerate complex tangents exactly along a single knot (of any given type), would be a worthy result. It is a problem that we would like to see fully resolved in the future. Possible generalizations of our results to the wider class of 3-manifolds are also of interest. 
\par\ \par\


\begin{thebibliography}{13}
\bibitem{[1]} P. Ahern and W. Rudin, ``Totally Real Embeddings of $S^3$ in $\mathbb{C}^3$,'' Proceedings of AMS , Vol. 94, No.5 (1985).

\bibitem{[2]} A. V. Domrin, "A description of characteristic classes of real submanifolds in complex manifolds via RC-singularities," Izv. Ross. Akad. Nauk Ser. Mat. (English transl.) (5) 59 (1995), 899918. MR 1360632 (97c:32020), Zbl 0877.57011

\bibitem{[3]} A. Elgindi, "On the Topological Structure of Complex Tangencies to Embeddings of $S^3$ into $\mathbb{C}^3$,"
New York J. of Math. Vol. 18 (2012), 295-313.

\bibitem{[4]} A. Elgindi, "On the Bishop invariants of embeddings
of $S^3$ into $\mathbb{C}^3$"
New York J. of Math. Vol. 20 (2014), 275-292.

\bibitem{[5]} A. Elgindi, "A Topological Obstruction to the Removal of a Degenerate Complex Tangent and Some Related Homotopy and Homology Groups," accepted for publication with the International Journal of Mathematics.


\bibitem{[6]} F. Forstneric, "On totally real embeddings into $\mathbb{C}^n$," Expositiones Mathematicae, 4 (1986), pp. 243255.

\bibitem{[7]} M. Gromov, ``Partial Differential Relations," Springer-Verlag (1986).

\bibitem{[8]} A. Hatcher, ``Algebraic Topology,'' Cambridge University Press (2002).

\bibitem{[9]} H.F. Lai, "Characteristic classes of real manifolds immersed in complex manifolds," Trans. AMS 172 (1972), 133.

\bibitem{[10]} J. Lee, ``Introduction to Smooth Manifolds,'' Springer Graduate Texts in Mathematics, Vol. 218 (2012).

\bibitem{[11]} M. Slapar, "Cancelling complex points in codimension two," Bull. Aust. Math. Soc. 88 (2013), no. 1, 6469.

\bibitem{[12]} S. Webster, ``The Euler and Pontrjagin numbers of an n-manifold in $\mathbb{C}^n$," Comm. Math. Helve. Vol. 60 (1985).




\end{thebibliography}
\end{document}